\documentclass[12pt]{article2}
\usepackage{amsmath}
\usepackage{amssymb}
\usepackage{amsthm}
\numberwithin{equation}{section}
\newtheorem{theorem}{Theorem}[section]

\newtheorem{Definition}[theorem]{Definition}

\newtheorem{Remark}[theorem]{Remark}
\newenvironment{remark}{\begin{Remark}\rm}{\end{Remark}}
\newtheorem{Example}[theorem]{Example}

\newcommand\beq{\begin{equation}}
\newcommand\eeq{\end{equation}}
\newcommand\bea{\begin{eqnarray}}
\newcommand\eea{\end{eqnarray}}
\newcommand\nonu{\nonumber\\}
\newcommand\sLP{\\[\smallskipamount]}
\newcommand\bLP{\\[\bigskipamount]}
\newcommand\PP{\par}
\newcommand\mPP{\\[\medskipamount]\indent}
\newcommand\bPP{\\[\bigskipamount]\indent}

\newcommand\RR{\mathbb{R}}
\newcommand\ZZ{\mathbb{Z}}
\newcommand\Zplus{\ZZ_{\ge 0}}
\newcommand\FSE{{\cal E}}
\newcommand\al\alpha
\newcommand\be\beta
\newcommand\ga\gamma
\newcommand\de\delta
\newcommand\tha\theta
\newcommand\la\lambda
\newcommand\si\sigma
\newcommand\Ga\Gamma
\newcommand\half{\frac12}
\newcommand\thalf{\textstyle\half}
\newcommand\quart{\frac14}
\newcommand\iy\infty
\newcommand\pa\partial
\renewcommand\Re{{\rm Re\,}}
\newcommand{\hyp}[5]{\,\mbox{}_{#1}F_{#2}\left(
  \genfrac{}{}{0pt}{}{#3}{#4};#5\right)}
\newcommand{\qhyp}[5]{\,\mbox{}_{#1}\phi_{#2}\left(
  \genfrac{}{}{0pt}{}{#3}{#4};#5\right)}
\newcommand\LHS{left-hand side}
\newcommand\RHS{right-hand side}
\newcommand\wrt{with respect to }
\newcommand\dstyle{\displaystyle}
\newcommand{\prat}[4]{p_{#1}^{(#2,#3)}(#4;q,{\rm rat})}

\begin{document}

\title{A second addition formula for continuous\\
$q$-ul\-tra\-spher\-i\-cal polynomials
\footnote{to appear in {\em Theory and Applications of Special Functions.
A Volume Dedicated to Mizan Rahman},\goodbreak
M.~E.~H. Ismail and E.~Koelink (eds.),
Developments in Mathematics, Kluwer}}
\author{Tom H. Koornwinder\\[\bigskipamount]
{\em Dedicated to Mizan Rahman}}

\date{June 26, 2003}

\maketitle

\begin{abstract}
This paper provides the details of Remark 5.4 in the author's paper
``Askey-Wilson polynomials as zonal spherical functions on the $SU(2)$
quantum group'', SIAM J.\ Math.\ Anal.\ 24 (1993), 795--813. 
In formula (5.9) of the 1993 paper a two-parameter class of
Askey-Wilson polynomials was expanded as a finite Fourier series with a
product of two ${}_3\phi_2$'s as Fourier coefficients.
The proof given there used
the quantum group interpretation. Here this
identity will be generalized to a
3-parameter class of Askey-Wilson polynomials being
expanded in terms of continuous $q$-ultraspherical polynomials with a
product of two ${}_2\phi_2$'s as coefficients,
and an analytic proof will be given for it.
Then Gegenbauer's addition formula for ultraspherical
polynomials and Rahman's addition formula for $q$-Bessel
functions will be obtained as limit cases.
This $q$-analogue of Gegenbauer's addition
formula is quite different from the addition formula for continuous
$q$-ultraspherical polynomials obtained by Rahman and Verma in 1986.
Furthermore, the functions occurring as factors
in the expansion coefficients will be interpreted
as a special case of a system of biorthogonal rational functions with
respect to the Askey-Roy $q$-beta measure.
A degenerate case of this biorthogonality are Pastro's biorthogonal
polynomials associated with the Stieltjes-Wigert polynomials.
\end{abstract}
\section{Introduction} \label{sec:1}
%until 1.22
%
Rahman and Verma \cite{RahV86} obtained the following {\em addition formula
for continuous $q$-ultraspherical polynomials\/}:
\bea
&&p_n(\cos\tha;a,aq^\half,-a,-aq^\half\mid q)=
\sum_{k=0}^n
\frac{(q;q)_n\,(a^4q^n,a^4q^{-1},a^2q^\half,-a^2q^\half,-a^2;q)_k\,a^{n-k}}
{(q;q)_k\,(q;q)_{n-k}\,(a^4 q^{-1};q)_{2k}\,
(a^2q^\half,-a^2q^\half,-a^2;q)_n}\,
\nonu
&&\hskip 3cm\times
p_{n-k}(\cos\phi;aq^{\half k},aq^{\half(k+1)},-aq^{\half k},-aq^{\half(k+1)}
\mid q)
\nonu
&&\hskip 3cm\times
p_{n-k}(\cos\psi;aq^{\half k},aq^{\half(k+1)},-aq^{\half k},-aq^{\half(k+1)}
\mid q)
\nonu
&&\hskip 5cm\times
p_k(\cos\tha;ae^{i(\phi+\psi)},ae^{-i(\phi+\psi)},ae^{i(\phi-\psi)},
ae^{i(\psi-\phi)}\mid q).
\label{eq:1.01}
\eea
The formula is here written in the form given in
\cite[Exercise 8.11]{GRah90}.
Use \cite{GRah90} also for notation of ($q$-)hypergeometric functions
and ($q$-)shifted factorials. Throughout it is
supposed that $0<q<1$.
\PP
Formula \eqref{eq:1.01} is given in terms of
{\em Askey-Wilson polynomials\/} (see \cite{AW85} or \cite[\S7.5]{GRah90}):
\beq
p_n(\cos\tha;a,b,c,d\mid q):=
a^{-n}\,(ab,ac,ad;q)_n\,r_n(\cos\tha;a,b,c,d\mid q)\quad
(n\in\Zplus)
\label{eq:1.02}
\eeq
(symmetric in $a,b,c,d$),
where
\beq
r_n(\cos\tha;a,b,c,d\mid q):=
\qhyp43
{q^{-n},abcdq^{n-1},ae^{i\tha},ae^{-i\tha}}{ab,ac,ad}{q,q}.
\label{eq:1.03}
\eeq
The {\em continuous $q$-ultraspherical polynomials\/} are the special case
$b=aq^\half, c=-a, d=-aq^\half$ of the Askey-Wilson polynomials, often
notated as follows (see \cite[(7.4.14)]{GRah90}):
\beq
C_n(x;a^2\mid q)=
\frac{(a^4;q)_n}{(q;q)_n a^n}\,
r_n(x;a,aq^\half,-a,-aq^\half\mid q)
=\frac{(a^2;q)_n}{(q,a^4q^n;q)_n}\,
p_n(x;a,aq^\half,-a,-aq^\half\mid q).
\label{eq:1.22}
\eeq
A further specialization to $a=\tfrac14$, i.e.,
$(a,b,c,d)=(\tfrac14,\tfrac34,-\tfrac14,-\tfrac34)$,
yields the {\em continuous $q$-Legendre polynomials}.
For this case $a=\tfrac14$ Koelink was able to give two different poofs of
the addition formula \eqref{eq:1.01} from a quantum group interpretation
on $SU_q(2)$, see \cite{Koel94} and \cite{Koel97}.
\PP
If $a$ is replaced by $q^{\half\la}$ in \eqref{eq:1.01} and the limit
is taken for $q\uparrow 1$, then a version of the
{\em addition formula for ultraspherical polynomials\/} is obtained:
\bea
&&C_n^\la(\cos\tha)=
\sum_{k=0}^n\frac{2^{2k}\,(2\la+2k-1)\,(n-k)!\,(\la)_k^2}
{(2\la-1)_{n+k+1}}\nonu
&&\qquad\quad\times
(\sin\phi)^k\,C_{n-k}^{\la+k}(\cos\phi)\,
(\sin\psi)^k\,C_{n-k}^{\la+k}(\cos\psi)\,
C_k^{\la-\half}\left(\frac{\cos\tha-\cos\phi\cos\psi}{\sin\phi\sin\psi}
\right).
\label{eq:1.04}
\eea
Here {\sl ultraspherical polynomials} are defined by
\beq
C_n^\la(x):=\frac{(2\la)_n}{n!}\,
\hyp21{-n,n+2\la}{\la+\thalf}{\thalf(1-x)}.
\label{eq:1.07}
\eeq
By elementary substitution the addition formula \eqref{eq:1.04}
transforms into the
familiar addition formula for ultraspherical polynomials:
\bea
&&C_n^\la(\cos\phi\cos\psi+\sin\phi\sin\psi\cos\tha)=
\sum_{k=0}^n
\frac{2^{2k}\,(2\la+2k-1)\,(n-k)!\,(\la)_k^2}{(2\la-1)_{n+k+1}}
\nonu
&&\qquad\qquad\qquad\qquad\times
(\sin\phi)^k\,C_{n-k}^{\la+k}(\cos\phi)\,
(\sin\psi)^k\,C_{n-k}^{\la+k}(\cos\psi)\,
C_k^{\la-1/2}(\cos\tha),
\label{eq:1.05}
\eea
see \cite[10.9(34)]{EMagOT53}, but watch out for the misprint $2^m$
which should be $2^{2m}$; see also the references given in
\cite[Lecture 4]{A75}.
For the removable singularity at $\la=\thalf$ in \eqref{eq:1.05}
observe that
\beq
\lim_{\la\to \half}\,\frac k{2\la-1}\,C_k^{\la-\half}(\cos\tha)=
\cos(k\tha)\quad(k>0);
\qquad
C_0^{\la-\half}(\cos\tha)=1.
\label{eq:1.08}
\eeq
\PP
An elementary transformation comparable to the passage from
\eqref{eq:1.04} to \eqref{eq:1.05} cannot be performed on the $q$-level.
It is a generally observed phenomenon in $q$-theory that, for a classical
formula involving a parameter dependent function with an
argument transformed by a parameter dependent
transformation, a possible $q$-analogue
has the transformation parameters occurring
in the function parameters. Compare for instance the
$p_k$ factor in \eqref{eq:1.01} with the $C_k^{\la-\half}$ factor
in~\eqref{eq:1.04}.
\PP
In \cite[(5.9)]{Koo93} I obtained the following formula as a spin-off of the
interpretation of certain Askey-Wilson polynomials as zonal spherical
functions on the quantum group $SU_q(2)$:
\bea
&&p_n(\cos\tha;-q^{\half(\si+\tau+1)},-q^{\half(-\si-\tau+1)},
q^{\half(\si-\tau+1)},q^{\half(-\si+\tau+1)}\mid q)\nonu
&&\qquad=
\sum_{k=-n}^n
\frac{(q^{n+1};q)_n (q;q)_{2n}}{(q;q)_{n+k} (q;q)_{n-k}}\,
q^{\half(n-k)(n-k+\si+\tau)}
\nonu
&&\qquad\qquad\quad\times
\qhyp32{q^{-n+k},q^{-n},q^{-n-\si}}{q^{-2n},0}{q,q}\,
\qhyp32{q^{-n+k},q^{-n},q^{-n-\tau}}{q^{-2n},0}{q,q}\,
e^{ik\tha}.
\label{eq:1.06}
\eea
Here the summand, with $e^{ik\tha}$ omitted,
is invariant under the transformation $k\to -k$,
as can be seen by twofold application of \cite[(3.2.3)]{GRah90}.
The ${}_3\phi_2$'s were viewed in \cite{Koo93} as dual $q$-Krawtchouk
polynomials (see their definition in \cite[\S3.17]{{KoekSw98}}),
but here we will
prefer to consider them as certain (unusual) $q$-analogues of
ultraspherical polynomials, for which an expression in terms of a ${}_2\phi_2$
is more suitable.
For this purpose apply a
${}_3\phi_2 \to {}_2\phi_2$ transformation obtained from
formulas (1.5.4) and (III.7) in \cite{GRah90}. Then, after the substitutions
$q^{\half\si}\to is^{-1}$, $q^{\half\tau}\to it$ in \eqref{eq:1.06},
we can write \eqref{eq:1.06} equivalently as
\bea
&&\hskip -0.5cm
p_n(\cos\tha;st^{-1}q^\half,s^{-1}tq^\half,-stq^\half,
-s^{-1}t^{-1}q^\half\mid q)
\nonu
&&\hskip -0.5cm
=(-1)^n q^{\half n^2}\,(q;q)_n
\sum_{k=0}^n q^{\half k}\,
\frac{(q^{n+1};q)_k\,(q^{-n};q)_k}{(q;q)_k^2}\,
(st)^{k-n}\,(-q^{-k+1}s^2,-q^{-k+1}t^{-2};q)_{n-k}
\nonu
&&\hskip -0.5cm\times
\qhyp22{q^{-n+k},q^{n+k+1}}{q^{k+1},-q^{k+1}s^2}{q,-qs^2}
\qhyp22{q^{-n+k},q^{n+k+1}}{q^{k+1},-q^{k+1}t^{-2}2}{q,-qt^{-2}}
(1+\de_{k,0})\cos(k\tha).\quad
\label{eq:1.21}
\eea
\PP
After the substitutions
$q^\si=\tan^2\thalf\phi$,
$q^\tau=\tan^2\thalf\psi$ in \eqref{eq:1.06},
the limit for $q\uparrow1$ becomes the case
$\la=\thalf$ of the addition formula \eqref{eq:1.05} (combined with
\eqref{eq:1.08}),
i.e., the addition formula for Legendre polynomials $P_n(x)=C_n^\half(x)$.
\PP
The first main result of this paper is the following addition formula for
continuous $q$-ultraspherical polynomials.
Thus I will
finally fulfill my promise of bringing out
``reference [9]'' of my paper \cite{Koo93}, which reference was mentioned
there as being in preparation. 
\begin{theorem} \label{th:1.12}
We have, in notation \eqref{eq:1.03},
\bea
&&\hskip -0.5cm
r_n(\cos\tha;ast^{-1}q^\half,ats^{-1}q^\half,-astq^\half,
-as^{-1}t^{-1}q^\half\mid q)
=(-1)^n a^{2n} q^{\half n(n+1)} \,(1+a^2)
\nonu
&&\hskip -0.5cm
\times\sum_{k=0}^n 
\frac{(1-a^2q^k)\,(q^{-n},a^4q,a^4q^{n+1};q)_k\,
(a^{-2}st^{-1}q^\half)^k}
{(1-a^4q^k)\,(q;q)_k\,(a^2q;q)_k^2\,(-a^2s^2q;q)_k\,(-a^2t^{-2}q;q)_k}\,
\qhyp22{q^{-n+k},a^4q^{n+k+1}}
{a^2q^{k+1},-a^2s^2 q^{k+1}}{q,-s^2q}
\nonu
&&\qquad\qquad\times
\qhyp22{q^{-n+k},a^4q^{n+k+1}}
{a^2q^{k+1},-a^2t^{-2}q^{k+1}}{q,-t^{-2}q}
r_k(\cos\tha;a,-a,a q^\half,-a q^\half\mid q),
\label{eq:1.09}
\eea
or, in notation \eqref{eq:1.02},
\bea
&&\hskip -0.5cm
p_n(\cos\tha;ast^{-1}q^\half,ats^{-1}q^\half,-astq^\half,
-as^{-1}t^{-1}q^\half\mid q)
=(-1)^nq^{\half n^2}\sum_{k=0}^n a^{n-k}\,(a^2q^{k+1};q)_{n-k}
\nonu
&&\hskip -0.5cm
\times
q^{\half k}\,
\frac{(q^{-n},a^4q^{n+1};q)_k}
{(q,a^4q^k;q)_k}\,
\frac{(-a^2s^2q^{k+1},-a^2t^{-2} q^{k+1};q)_{n-k}}{(st^{-1})^{n-k}}
\qhyp22{q^{-n+k},a^4q^{n+k+1}}
{a^2q^{k+1},-a^2s^2 q^{k+1}}{q,-s^2q}
\nonu
&&\qquad\qquad\times
\qhyp22{q^{-n+k},a^4q^{n+k+1}}
{a^2q^{k+1},-a^2t^{-2}q^{k+1}}{q,-t^{-2}q}
p_k(\cos\tha;a,-a,a q^\half,-a q^\half\mid q).
\label{eq:1.20}
\eea
\end{theorem}
\PP
For $a=1$ formula \eqref{eq:1.20} specializes
to formula \eqref{eq:1.21} (with usage of a Chebyshev case of the
Askey-Wilson polynomials,
see \cite[(4.21)]{AW85}). Its limit case for
$q\uparrow1$ (after the substitutions
$a=q^{\half\la-\quart}$, $s=\tan(\half\phi)$, $t=\tan(\half\psi)$)
is the addition formula \eqref{eq:1.05} for ultraspherical polynomials of
general order $\la$.
It is also possible to obtain Rahman's addition formula
\cite[(1.10)]{Rah88}
for $q$-Bessel functions as a formal limit case of
\eqref{eq:1.20}, see Remark \ref{th:2.23}.
For precise correspondence of \eqref{eq:1.09} with \eqref{eq:1.01}
one should replace $a$ by $aq^{-\quart}$ in \eqref{eq:1.09}.
\PP
The \LHS\ of \eqref{eq:1.20} is invariant under the symmetries
$(s,t)\to(t,s)$, $(s,t)\to (-s^{-1},t)$ and $(s,t)\to (s,-t^{-1})$.
These symmetries are also visible on the \RHS\ of \eqref{eq:1.20}
if we take into account that
\[
\frac{(-a^2s^2q^{k+1};q)_{n-k}}{s^{n-k}}\,\qhyp22{q^{-n+k},a^4q^{n+k+1}}
{a^2q^{k+1},-a^2s^2 q^{k+1}}{q,-s^2q}
\]
is invariant under the transformation $s\to s^{-1}$ up to the
factor $(-1)^{n-k}$ (by \eqref{eq:3.14}).
\PP
The proof of Theorem \ref{th:1.12}
(see details in \S\ref{sec:2})
is quite similar to the proof
of \eqref{eq:1.01} in \cite{RahV86}. We consider \eqref{eq:1.09} as a
connection formula which connects Askey-Wilson polynomials of different
order. There are certain choices of the orders of the Askey-Wilson
polynomials in a connection formula
\beq
r_n(\cos\tha;\al,\be,\ga,\de\mid q)=
\sum_{k=0}^n A_{n,k}\,r_k(\cos\tha;a,b,c,d\mid q)
\label{eq:1.10}
\eeq
for which the connection coefficients $A_{n,k}$ factorize.
Formula \eqref{eq:1.01} is one example; formula  \eqref{eq:1.09}
is another example.
\PP
At the end of \S2 a degenerate addition formula for
continuous $q$-ultraspherical polynomials will be given as a limit case
of the addition formula \eqref{eq:1.09}.
As a further limit case we obtain a degenerate addition formula
for $q$-Bessel functions.
\PP
The ${}_2\phi_2$ factors on the \RHS\ of \eqref{eq:1.09} tend, 
after the mentioned
substitutions and as $q\uparrow1$,
to ul\-tra\-spher\-i\-cal polynomials $C_{n-k}^{\la+k}$ of argument
$\cos\phi$ resp.\ $\cos\psi$, so we might expect that these ${}_2\phi_2$'s
also
satisfy some kind of (bi)orthogonality relations.
This is indeed the case $\al=\be$ 
of Theorem \ref{th:1.13} below. See its proof in \S\ref{sec:3}.
\begin{theorem} \label{th:1.13}
Let $\al,\be>-1$.
Define a system of rational functions in $t$ by
\beq
\prat n\al\be t:=
\qhyp22{q^{-n},q^{n+\al+\be+1}}
{q^{\al+1},-tq^{\be+1}}{q,-tq}
\qquad(n\in\Zplus),
\label{eq:1.14}
\eeq
and define, with additional parameter $c\in\RR$, a measure on $[0,\iy)$ by
\beq
d\mu_{\al,\be,c;q}(t):=
\frac{q^{-c^2}\,\Ga_q(c)\,\Ga_q(1-c)\,\Ga_q(\al+\be+2)}
{\Ga(c)\,\Ga(1-c)\,\Ga_q(\al+1)\,\Ga_q(\be+1)}\,
\frac{t^{c-1}\,(-tq^{\be+1},-t^{-1}q^{\al+2};q)_\iy}
{(-tq^{-c},-t^{-1}q^{1+c};q)_\iy}\,dt,
\label{eq:1.15}
\eeq
where $\Ga_q$ is defined in \eqref{eq:3.02}.
Then
\beq
\int_0^\iy
\prat n\al\be t\,\prat m\be\al{t^{-1}q}\,
d\mu_{\al,\be,c;q}(t)=
\frac{(-1)^n\,(1-q^{n+\al+\be+1})\,(q;q)_n\,\de_{n,m}}
{q^{\half n(n-1)}\,(1-q^{2n+\al+\be+1})\,(q^{\al+\be+2};q)_n}\,.
\label{eq:1.16}
\eeq
\end{theorem}
\PP
The case $n=m=0$ of \eqref{eq:1.16}, i.e., the integration formula
$\int_0^\iy d\mu_{\al,\be,c;q}(t)=1$, is precisely the $q$-beta integral
of Askey and Roy \cite[(3.4)]{ARo86}, which extends a
$q$-beta integral of Ramanujan \cite[(19)]{Ram15}.
The Askey-Roy integral was
independently obtained by Gasper (see \cite{G84} and also
\cite{G87}) and by Thiruvenkatachar \& Venkatachaliengar
(see \cite[p.93]{A88}).
\PP
Note that the two ${}_2\phi_2$'s in \eqref{eq:1.09} (with $a=q^{\half\al}$),
can be rewritten in
terms of the function \eqref{eq:1.14}: as
$\prat{n-k}{\al+k}{\al+k}{s^2}$ and $\prat{n-k}{\al+k}{\al+k}{t^{-2}}$,
respectively.
\PP
The biorthogonality measure in \eqref{eq:1.16} is evidently not unique,
because of the parameter $c$.
Further illustration of the non-uniqueness of the measure for
these biorthogonality relation is provided by a $q$-integral variant of
\eqref{eq:1.16}. In order to state this, we need the following
definition of $q$-integral on $(0,\iy)$ ($f$ arbitrary function on $(0,\iy)$
for which the sum absolutely converges):
\beq
\int_0^{s\cdot\iy} f(t)\,d_qt:=
(1-q)\sum_{k=-\iy}^\iy sq^k\,f(sq^k)\qquad(s>0).
\label{eq:1.18}
\eeq
\begin{theorem} \label{th:1.17}
The functions defined by \eqref{eq:1.14} also satisfy
the biorthogonality relations
\bea
&&\frac{\Ga_q(\al+\be+2)}{\Ga_q(\al+1)\,\Ga_q(\be+1)}\,
\int_0^{s\cdot\iy}
\prat n\al\be t\,\prat m\be\al{t^{-1}q}\,
\frac{t^{-1}\,(-tq^{\be+1},-t^{-1}q^{\al+2};q)_\iy}
{(-t,-t^{-1}q;q)_\iy}\,d_qt
\nonu
&&\qquad\qquad\qquad\qquad=
\frac{(-1)^n\,(1-q^{n+\al+\be+1})\,(q;q)_n\,\de_{n,m}}
{q^{\half n(n-1)}\,(1-q^{2n+\al+\be+1})\,(q^{\al+\be+2};q)_n}
\qquad(\al,\be>-1,\;s>0).\quad
\label{eq:1.19}
\eea
\end{theorem}
\PP
The case $n=m=0$ of \eqref{eq:1.19} is a $q$-beta integral first given
by Gasper \cite{G87}, but it is essentially Ramanujan's ${}_1\psi_1$ sum;
see some further discussion in \S\ref{sec:3}.
The proof of Theorem \ref{th:1.17} is by a completely analogous argument
as I will give in \S\ref{sec:3} for the proof of Theorem \ref{th:1.13}.
\PP
The paper concludes in \S\ref{sec:4}
with some open questions and with some specializations of Theorem
\ref{th:1.13}. Pastro's \cite{P85} biorthogonal polynomials associated
with the Stieltjes-Wigert polynomials occur as a special case.
\bLP
{\em Acknowlegement}\quad
I did the work communicated here essentially already in the beginning
of the nineties. During that time
Mizan Rahman sent me very useful hints concerning
the material which is now in \S2, while Ren\'e Swarttouw
carefully checked
(and corrected) my computations.
One of the referees made some very good suggestions,
which resulted, among others, into Remarks \ref{th:2.22} and
\ref{th:2.25}.
Finally I thank Erik Koelink
for stimulating me to publish this work after such a long time.
\section{Proof of the new addition formula} \label{sec:2}
%until 2.28
%
In this section I prove the second addition formula for
continuous $q$-ultraspherical polynomials, stated in
Theorem \ref{th:1.12}.
Let us first consider the general connection formula \eqref{eq:1.10}.
We can split up this connection into three successive connections of more
simple nature:
\bea
r_n(\cos\tha;\al,\be,\ga,\de\mid q)&=&
\sum_{j=0}^n c_{n,j}\,(\al e^{i\tha},\al e^{-i\tha};q)_j\,,
\label{eq:2.01}\\
(\al e^{i\tha},\al e^{-i\tha};q)_j&=&
\sum_{l=0}^j d_{j,l}\,(a e^{i\tha},a e^{-i\tha};q)_l\,,
\label{eq:2.02}\\
(ae^{i\tha},ae^{-i\tha};q)_l&=&
\sum_{k=0}^l e_{l,k}\,r_k(\cos\tha;a,b,c,d\mid q).
\label{eq:2.03}
\eea
Then
\beq
A_{n,k}=
\sum_{j=0}^{n-k}\sum_{l=0}^j c_{n,j+k}\,d_{j+k,l+k}\,e_{l+k,k}
\label{eq:2.04}
\eeq
and
\beq
A_{n,k}=\sum_{m=0}^{n-k}\sum_{i=0}^m c_{n,n-i}\,d_{n-i,n-m}\,e_{n-m,k}.
\label{eq:2.05}
\eeq
The coefficients $c_{n,j}$, $d_{j,l}$ and $e_{l,k}$ can be explicitly given as
\bea
c_{n,j}&=&\frac{(q^{-n},\al\be\ga\de q^{n-1};q)_j\,q^j}
{(\al\be,\al\ga,\al\de,q;q)_j}\,,
\label{eq:2.06}
\\
d_{j,l}&=&\frac{(q,a\al;q)_j\,(a^{-1}\al;q)_{j-l}\,\al^l}
{(q,a\al;q)_l\,(q;q)_{j-l}\,a^l}\,,
\label{eq:2.07}
\\
e_{l,k}&=&\frac{(abcd;q)_{2k}\,(q^{-l};q)_k\,(ab,ac,ad;q)_l\,q^{lk}}
{(abcdq^{k-1},q;q)_k\,(abcd;q)_{k+l}}\,.
\label{eq:2.08}
\eea
Here \eqref{eq:2.06} follows from \eqref{eq:1.03},
formula \eqref{eq:2.07} can be obtained by rewriting the
$q$-Saalsch\"utz formula \cite[(1.7.2)]{GRah90},
and \eqref{eq:2.08} follows from \cite[(2.6), (2.5)]{AW85}.

It turns out that in the two double sums \eqref{eq:2.04}, \eqref{eq:2.05}
of $A_{n,k}$ the inner sum can be written as a balanced ${}_4\phi_3$ of
argument $q$:
\bea
&&A_{n,k}=
(-1)^k\,q^{\half k(k+1)}\,
\frac{(q^{-n},\al\be\ga\de q^{n-1},ab,ac,ad;q)_k\,\al^k}
{(q,abcdq^{k-1},\al\be,\al\ga,\al\de;q)_k\,a^k}
\nonu
&&\quad\times\sum_{j=0}^{n-k}
\frac{(q^{-n+k},\al\be\ga\de q^{n+k-1},\al aq^k,\al a^{-1};q)_j\,q^j}
{(\al\be q^k,\al\ga q^k,\al\de q^k,q;q)_j}\,
\qhyp43{q^{-j},abq^k,acq^k, adq^k}
{a\al q^k,a\al^{-1}q^{1-j},abcdq^{2k}}{q,q}\qquad
\label{eq:2.09}
\eea
and
\bea
&&\hskip-1.2cm A_{n,k}=
(-1)^n\,q^{-\half n(n-2k-1)}\,
\frac
{(1-abcdq^{2k-1})\,(q^{-n};q)_k\,(\al\be\ga\de q^{n-1},ab,ac,ad;q)_n\,\al^n}
{(1-abcdq^{k-1})\,(q;q)_k\,(abcdq^k,\al\be,\al\ga,\al\de;q)_n\,a^n}
\nonu
&&
\qquad\quad\times\sum_{m=0}^{n-k}
\frac{(q^{-n+k},(abcd)^{-1}q^{1-n-k},\al a^{-1},(\al a)^{-1}q^{1-n};q)_m\,q^m}
{(q,(ab)^{-1}q^{1-n},(ac)^{-1}q^{1-n},(ad)^{-1}q^{1-n};q)_m}
\nonu
&&\qquad\qquad\qquad\times
\qhyp43
{q^{-m},(\al\be)^{-1}q^{-n+1},(\al\ga)^{-1}q^{-n+1},(\al\de)^{-1}q^{-n+1}}
{(\al\be\ga\de)^{-1}q^{-2n+2},(\al a)^{-1}q^{-n+1},a\al^{-1}q^{-m+1}}{q,q}\,.
\label{eq:2.10}
\eea
The sums in \eqref{eq:2.09} and \eqref{eq:2.10} can be compared with
the {\em Bateman type product formula\/} \cite[(8.4.7)]{GRah90}:
\bea
&&r_n(\cos\phi;a,aq^\half,-a,-aq^\half\mid q)\,
r_n(\cos\psi;a,aq^\half,-a,-aq^\half\mid q)
=
\frac{1+a^2q^n}{1+a^2}\,(-1)^n\,q^{-\half n}
\nonu
&&\qquad\qquad\times\sum_{m=0}^n
\frac
{(q^{-n},a^4q^n,-e^{i\psi-i\phi}q^\half,-a^2e^{i\phi-i\psi}q^\half;q)_m\,q^m}
{(q,a^2q^\half,-a^2q^\half,-a^2q;q)_m}
\nonu
&&\qquad\qquad\qquad\qquad\qquad\qquad\times
\qhyp43{q^{-m},a^2,a^2e^{2i\phi},a^2e^{-2i\psi}}
{a^4,-a^2e^{i\phi-i\psi}q^\half,-e^{i\phi-i\psi}q^{-m+\half}}{q,q}.
\label{eq:2.11}
\eea
The sum in \eqref{eq:2.09} can be matched with the \RHS\ of
\eqref{eq:2.11} precisely for those values of the parameters
$a,b,c,d,\al,\be,\ga,\de$ in \eqref{eq:1.10}
which occur in the Rahman-Verma
addition formula \eqref{eq:1.01},
i.e., for $ab=cd$, $\al^2=cdq$, $\al=\be q^\half=-\ga q^\half=-\de$.
In fact, this will prove \eqref{eq:1.01}. The proof in
\cite{RahV86} is essentially along these lines.
\PP
Next we see that the sum in \eqref{eq:2.10} can be matched with
the \RHS\ of \eqref{eq:2.11} precisely for those values of the parameters
$a,b,c,d,\al,\be,\ga,\de$ in \eqref{eq:1.10} when
$\al\be=\ga\de=a^2q$,
$b=-c=aq^\half=-dq^\half$.
Thus put $b=aq^\half$, $c=-aq^\half$, $d=-a$,
$\al=ast^{-1}q^\half$,
$\be=ats^{-1}q^\half$,
$\ga=-astq^\half$,
$\de=-as^{-1}t^{-1}q^\half$ in \eqref{eq:1.10} and \eqref{eq:2.10}.
Then these two formulas specialize to:
\beq
r_n(\cos\tha;ast^{-1}q^\half,ats^{-1}q^\half,-astq^\half,-as^{-1}t^{-1}q^\half
\mid q)
=\sum_{k=0}^n A_{n,k}\,
r_k(\cos\tha;a,-a,a q^\half,-a q^\half\mid q),
\label{eq:2.12}
\eeq
\bea
&&A_{n,k}=
(-1)^n\,s^n\,t^{-n}\,q^{\half n(-n+2k+2)}\,
\frac{(1+a^2)\,(1-a^4q^{2k})\,(q^{-n},a^4q;q)_k\,(a^4q^{n+1};q)_n^2}
{(1+a^2q^n)\,(1-a^4q^k)\,(q,a^4q^{n+1};q)_k\,(a^2q;q)_n^2}
\nonu
&&\qquad\times
\frac1{(-s^2a^2q,-t^{-2}a^2q;q)_n}\,
\sum_{m=0}^{n-k}
\frac{(q^{-n+k},a^{-4}q^{-n-k},st^{-1}q^\half,
a^{-2}ts^{-1}q^{-n+\half};q)_m\,q^m}
{(q,a^{-2}q^{\half-n},-a^{-2}q^{\half-n},-a^{-2}q^{1-n};q)_m}
\nonu
&&\qquad\qquad\qquad\qquad\times
\qhyp43
{q^{-m},a^{-2}q^{-n},-a^{-2}s^{-2}q^{-n},-a^{-2}t^2q^{-n}}
{a^{-4}q^{-2n},a^{-2}ts^{-1}q^{-n+\half},ts^{-1}q^{-m+\half}}{q,q}.
\label{eq:2.13}
\eea
If we next make the successive substitutions
$n\to n-k$, $a\to a^{-1}q^{-\half n}$, $e^{i\phi}\to is^{-1}$,
$e^{i\psi}\to -it^{-1}$ in \eqref{eq:2.11} and compare with
\eqref{eq:2.13} then we can write \eqref{eq:2.12} as follows:
\bea
&&r_n(\cos\tha;
ast^{-1}q^\half,ats^{-1}q^\half,-astq^\half,-as^{-1}t^{-1}q^\half\mid q)
\nonu
&&\qquad=\frac{s^n t^{-n} q^{-\half n(n-1)} (a^4q^{n+1};q)_n^2}
{(a^2q;q)_n^2\,(-s^2a^2q,-t^{-2}a^2q;q)_n}
\sum_{k=0}^n (-1)^k q^{k(n+\half)}\, \frac{1-a^2q^k}{1-a^2}\,
\frac{(q^{-n},a^4;q)_k}{(q,a^4q^{n+1};q)_k}\,
\nonu
&&\qquad\qquad\times
r_{n-k}(\tfrac{s-s^{-1}}{2i};
a^{-1}q^{-\half n},a^{-1}q^{-\half n+\half},-a^{-1}q^{-\half n},
-a^{-1}q^{-\half n+\half}\mid q)
\nonu
&&\qquad\qquad\times
r_{n-k}(\tfrac{t^{-1}-t}{2i};\,
a^{-1}q^{-\half n},a^{-1}q^{-\half n+\half},-a^{-1}q^{-\half n},
-a^{-1}q^{-\half n+\half}\mid q)\,
\nonu
&&\qquad\qquad\qquad\qquad\qquad\qquad\times
r_k(\cos\tha;a,-a,aq^\half,-aq^\half\mid q).
\label{eq:2.17}
\eea
\PP
In order to make the two $r_{n-k}$ factors on the \RHS\ above into
closer $q$-analogues of the two $C_{n-k}^{\la+k}$ factors on the \RHS\
of \eqref{eq:1.05}, we will use the following string of identities:
\bea
&&r_n(\cos\tha;a,aq^\half,-a,-aq^\half\mid q)
=\qhyp43{q^{-n},a^4 q^{n-1},ae^{i\tha},ae^{-i\tha}}
{a^2q^\half,-a^2,-a^2q^\half}{q,q}
\nonu
&&\qquad\qquad
=\frac{(a^2;q)_n\,a^n\,e^{in\tha}}{(a^4;q)_n}\,
\qhyp21{q^{-n},a^2}{a^{-2}q^{1-n}}{q,a^{-2}e^{-2i\tha}q}
\nonu
&&\qquad\qquad
=\frac{(a^2,a^{-2}e^{-2i\tha}q^{1-n};q)_n\,a^n\,e^{in\tha}}
{(a^4;q)_n}\,
\qhyp22{q^{-n},a^{-4}q^{1-n}}
{a^{-2}q^{1-n},a^{-2}e^{-2i\tha}q^{1-n}}{q,e^{-2i\tha}q}.\qquad
\label{eq:2.14}
\eea
For the proof use successively
(7.4.14), (7.4.2) and (1.5.4) in \cite{GRah90}.
\bLP
{\bf Proof of Theorem \ref{th:1.12}}
\sLP
This follows from \eqref{eq:2.17} by twofold substitution
of \eqref{eq:2.14}. Here replace $n$ by $n-k$ and $a$ by
$a^{-1}q^{-\half n}$ in \eqref{eq:2.17}, and replace $e^{i\tha}$ by 
$is^{-1}$ for the first substitution and by $it$ for the second
substitution.\qed
\begin{remark}
\label{th:2.22}
Let the divided difference operator $D_q$ acting on a function $F$ of
argument $e^{i\tha}$ be given by
\[
D_qF(e^{i\tha}):=\frac{\de_qF(e^{i\tha})}{\de_q\cos\tha},\quad
\de_qF(e^{i\tha}):=F(q^{\half}e^{i\tha})-F(q^{-\half}e^{i\tha}).
\]
Then $D_q$, acting on both sides of the addition formula \eqref{eq:1.20},
sends this to the same formula with $n$ replaced by $n-1$ and $a$ replaced
by $q^\half a$
(apply \cite[(3.1.9)]{KoekSw98}).
Thus, if formula \eqref{eq:1.20} is already known for $a=1$
(i.e., if formula \eqref{eq:1.21} is known) then the
procedure just sketched yields this formula for $a=q^{\half j}$ for all
$j\in\Zplus$, i.e., for infinitely many disctinct values of $a$.
Since, for fixed $n$, both sides of \eqref{eq:1.20} are rational in $a$,
formula \eqref{eq:1.20} will then be valid for general $a$.
Since formula \eqref{eq:1.06}, equivalent to \eqref{eq:1.21},
can be obtained by a quantum group
interpretation, we can say that it is possible to prove
formula \eqref{eq:1.20} by
arguments in
a quantum group setting, followed by minor analytic, but not very
computational reasoning.
\end{remark}
\noindent
{\bf Proof that \eqref{eq:1.09} has limit \eqref{eq:1.05}
as $q\uparrow1$}
\sLP
(after the substitutions
$a=q^{\half\la-\quart}$, $s=\tan(\half\phi)$, $t=\tan(\half\psi)$).
\sLP
The limits for $q\uparrow1$ of the factors
$r_n(\cos\tha)$, $r_k(\cos\tha)$, ${}_2\phi_2(-s^2q)$ and
${}_2\phi_2(-t^{-2}q)$ in \eqref{eq:1.09}, after the above substitutions and
after substitution of \eqref{eq:1.03} for the $r_n$ and $r_k$ factors
yields respectively:
\begin{eqnarray*}
{}_2F_1\bigl(-n,n+2\la;\la+\thalf;\thalf-\thalf
(\cos\phi\cos\psi+\sin\phi\sin\psi\cos\tha)\bigr),\;
{}_2F_1(-k,k+2\la-1;\la;\thalf-\thalf\cos\tha),
\\
{}_2F_1(-n+k,2\la+n+k;\la+k+\thalf;\thalf-\thalf\cos\phi),\;
{}_2F_1(-n+k,2\la+n+k;\la+k+\thalf;\thalf+\thalf\cos\psi).
\end{eqnarray*}
Express these ${}_2F_1$'s as ultraspherical polynomials by \eqref{eq:1.07}.
The limit of the coefficients on the \RHS\ of \eqref{eq:1.09}
(after the above substitutions) is also easily computed.\qed
\begin{remark}
\label{th:2.23}
Replace $s$ by $sq^{\half n}$ and $t$ by $tq^{-\half n}$ in
\eqref{eq:1.20} and let $n\to\iy$.
Then formally we obtain Rahman's addition formula for $q$-Bessel functions,
see \cite[(1.10)]{Rah88} (but watch out for the misprint $\Ga_q(\nu+1)$ which
should be $\Ga_q^2(\nu+1)$):
\bea
&&\qhyp21{-astq^\half e^{i\tha},-astq^\half e^{-i\tha}}
{a^2q}{q,-a^{-2}t^{-2}}=
\frac1{(-a^{-2}t^{-2};q)_\iy}\,
\sum_{k=0}^\iy
\frac{(-1)^kq^{\half k^2} a^{-k} s^k t^{-k}}
{(q,a^2q,a^4q^k;q)_k}
\nonu
&&\times
\qhyp01{-}{a^2q^{k+1}}{q,-s^2q^{k+1}}
\qhyp01{-}{a^2q^{k+1}}{q,-t^{-2}q^{k+1}}
p_k(\cos\tha;a,-a,aq^\half,-aq^\half\mid q).\qquad\quad
\label{eq:2.24}
\eea
According to \cite[(1.10)]{Rah88}, the further conditions $0<a<1$,
$0<s<t^{-1}$, $\tha\in\RR$ should be imposed here.
If we replace in \eqref{eq:2.24} $s$ by $(1-q)s$, $t$ by $(1-q)^{-1}t$,
$a$ by $q^{\half\al}$, and if we let $q\uparrow1$ then we formally
obtain the familiar Gegenbauer's
addition formula for Bessel functions $J_\al$,
see \cite[7.15(32)]{EMagOT53}.

In \eqref{eq:2.24} the \LHS\ is an {\em Askey-Wilson $q$-Bessel function}
(see \cite[\S2.3]{KoelSt01}), earlier studied under the name
{\em $q$-Bessel function on a $q$-quadratic grid} in
\cite{IMasSu99} and \cite{BuSu98}. The ${}_0\phi_1$'s on the \RHS\ of
\eqref{eq:2.24} are {\em Jackson's second $q$-Bessel functions},
usually written (see \cite[Exercise 1.24]{GRah90}) as
\beq
J_\al^{(2)}(x;q)=\frac{(q^{\al+1};q)_\iy}{(q;q)_\iy}\,(x/2)^\al\,
\qhyp01-{q^{\al+1}}{q,-\tfrac14 x^2 q^{\al+1}}.
\label{eq:2.28}
\eeq

Koelink \cite[(3.6.18)]{Koel91} gave an interpretation of the case $a=1$ of
\eqref{eq:2.24} on the quantum group of plane motions. This interpretation
is similar to the interpretation given
in \cite[(5.9)]{Koo93} for equation \eqref{eq:1.06}.
\end{remark}
\begin{remark}
\label{th:2.21}
If we multiply both sides of the addition formula \eqref{eq:1.09}
with $(-a^2t^{-2}q;q)_n$ and if we next let $t\to iaq^{\half n}$
in \eqref{eq:1.09} then we obtain a degenerate form of the addition
formula \eqref{eq:1.09}:
\bea
&&\frac{(-ise^{i\tha}q^{\half(1-n)},-ise^{-i\tha}q^{\half(1-n)};q)_n}
{(-a^2s^2q;q)_n}
=\sum_{k=0}^n(-ia^{-1}q^{\half(n+1)})^k\,
\frac{(q^{-n},a^4;q)_k}{(q,a^2;q)_k}\,
\frac{s^k}{(-a^2s^2q;q)_k}
\nonu
&&\quad\qquad\times
\qhyp22{q^{-n+k},a^4q^{n+k+1}}{a^2q^{k+1},-a^2s^2q^{k+1}}
{q,-s^2q}\,
r_k(\cos\tha;a,aq^\half,-aq^\half,-a\mid q).
\label{eq:2.20}
\eea
Integrated forms of \eqref{eq:1.09} and \eqref{eq:2.20} can be
obtained by integrating both sides of these formulas \wrt
the measure
$\left|\frac{(e^{2i\tha};q)_\iy}{(a^2 e^{2i\tha};q)_\iy}\right|^2d\tha$
on $[0,\pi]$, i.e., \wrt the orthogonality measure
for the continuous $q$-ultraspherical polynomials
$r_k(\cos\tha;a,q^\half a,-q^\half a,-a\mid q)$
(see \cite[\S7.4]{GRah90}).
This will yield a product formula and an integral representation,
respectively, for the functions $\prat n\al\al{s^2}$ (with $a=q^{\half\al}$).
\end{remark}
\begin{remark}
\label{th:2.25}
In \eqref{eq:2.20} replace $s$ by $sq^{\half n}$ and let $n\to\iy$.
Then we formally obtain a degenerate addition formula for $q$-Bessel
functions:
\bea
&&(-isq^\half e^{i\tha},-isq^\half e^{-i\tha};q)_\iy=
\sum_{k=0}^\iy i^k a^{-k} s^k q^{\half k^2}
\frac{(a^4;q)_k}{(q,a^2;q)_k}
\nonu
&&\qquad\qquad
\times\qhyp01{-}{a^2q^{k+1}}{q,-s^2q^{k+1}}
r_k(\cos\tha;a,aq^\half,-aq^\half,-a\mid q)
\label{eq:2.26}
\eea
If we replace in \eqref{eq:2.26} $a$ by $q^{\half\al}$, and if we substitute
\eqref{eq:2.28} and \eqref{eq:1.22}
then we can rewrite \eqref{eq:2.26} as
\bea
&&(-isq^\half e^{i\tha},-isq^\half e^{-i\tha};q)_\iy=
\frac{q^{\half\al^2}}{s^\al}\,
\frac{(q;q)_\iy}{(q^{\al+1};q)_\iy}
\sum_{k=0}^\iy i^k q^{\half k^2+\half k\al}\,\frac{1-q^{\al+k}}{1-q^\al}
\nonu
&&
\qquad\qquad\qquad\qquad\times
J_{\al+k}^{(2)}(2sq^{-\half\al};q)\,C_k(\cos\tha;q^\al\mid q).
\label{eq:2.27}
\eea
It is interesting to compare formula \eqref{eq:2.27} with
Ismail \& Zhang \cite[(3.32)]{IZha94}. They expand there a
generalized $q$-exponential function
$\FSE_q(z;-i,b/2)$ in terms of the $C_k(z;q^\al\mid q)$ and they obtain
almost the same expansion coefficients as in \eqref{eq:2.27}, including
Jackson's second $q$-Bessel functions, but they have a factor
$q^{\quart k^2}$, where \eqref{eq:2.27} has a factor
$q^{\half k^2+\half k\al}$.
\end{remark}
\section{Rational biorthogonal functions for the
Askey-Roy q-beta measure} \label{sec:3}
%until 3.15
%
The {\em Askey-Roy $q$-beta integral\/}
(see \cite[(3.4)]{ARo86}, \cite{G87},
\cite[pp. 92,93]{A88}
\cite[Exercise 6.17(ii)]{GRah90})
is as follows:
\bea
&&\int_0^\iy t^{c-1}\,
\frac{(-tq^{\be+1},-t^{-1}q^{\al+2};q)_\iy}
{(-tq^{-c},-t^{-1}q^{1+c};q)_\iy}\,dt=
\frac{q^{c^2}\,\Ga(c)\,\Ga(1-c)\,\Ga_q(\al+1)\,\Ga_q(\be+1)}
{\Ga_q(c)\,\Ga_q(1-c)\,\Ga_q(\al+\be+2)}\,,\qquad\qquad
\nonu
&&\qquad\qquad\qquad\qquad\qquad\qquad\qquad\qquad\qquad\qquad\qquad\qquad
\Re\al,\Re\be>-1,\;c\in\RR.
\label{eq:3.01}
\eea
Here the {\em $q$-gamma function\/} is defined by
\beq
\Ga_q(z):=\frac{(q;q)_\iy}{(q^z;q)_\iy}\,(1-q)^{1-z}.
\label{eq:3.02}
\eeq
The special case $c=\al+1$ of \eqref{eq:3.01} goes back (without proof) to
Ramanujan in Chapter 16 of his second notebook (see
\cite[p.29, Entry 14]{Be91}), and later in his paper
\cite[(19)]{Ram15}, with subsequent proof by
Hardy \cite{H15}.
\par
When we let $q\uparrow1$ in \eqref{eq:3.01} then we formally obtain the
{\em beta integral\/} on $(0,\iy)$:
\beq
\int_0^\iy t^\al\,(1+t)^{-\al-\be-2}\,dt=
\frac{\Ga(\al+1)\,\Ga(\be+1)}{\Ga(\al+\be+2)}\qquad
(\Re\al,\Re\be>-1).
\label{eq:3.03}
\eeq
When we move the orthogonality relations
\bea
&&\int_{-1}^1 P_n^{(\al,\be)}(x)\,P_m^{(\al,\be)}(x)\,(1-x)^\al\,(1+x)^\be\,dx
=\frac{2^{\al+\be+1}\,\Ga(n+\al+1)\,\Ga(n+\be+1)}
{(2n+\al+\be+1)\,n!\,\Ga(n+\al+\be+1)}\,\de_{n,m}
\nonu
&&\qquad\qquad\qquad\qquad\qquad\qquad\qquad\qquad\qquad\qquad
\qquad\qquad\qquad\qquad\qquad
(\al,\be>-1)\quad
\label{eq:3.04}
\eea
of the {\em Jacobi polynomials}
\beq
P_n^{(\al,\be)}(x):=
\frac{(\al+1)_n}{n!}\,
\hyp21{-n,n+\al+\be+1}{\al+1}{\thalf(1-x)}
\label{eq:3.05}
\eeq
(see \cite[\S10.8]{EMagOT53})
to $[0,\iy)$ by the substitution $x=(1-t)/(1+t)$,
then we obtain orthogonality relations
\bea
&&\int_0^\iy P_n^{(\al,\be)}\left(\frac{1-t}{1+t}\right)\,
P_m^{(\al,\be)}\left(\frac{1-t}{1+t}\right)\,
t^\al\,(1+t)^{-\al-\be-2}\,dt
\nonu
&&\qquad\qquad\qquad\qquad=
\frac{\Ga(n+\al+1)\,\Ga(n+\be+1)}
{(2n+\al+\be+1)\,n!\,\Ga(n+\al+\be+1)}\,\de_{n,m}\qquad
(\al,\be>-1).\quad
\label{eq:3.06}
\eea
Thus the rational functions $t\mapsto P_n^{(\al,\be)}\bigl((1-t)/(1+t)\bigr)$,
$n\in\Zplus$, are orthogonal \wrt the beta measure on $[
0,\iy)$ of which
the total mass is given in \eqref{eq:3.03}.
We would like to find $q$-analogues of these orthogonal rational functions
such that the orthogonality measure is the $q$-beta measure
in \eqref{eq:3.01}.
\PP
From \eqref{eq:1.15} we have
\beq
d\mu_{\al,\be,c;q}(t)=
\frac{w_{\al,\be,c;q}(t)\,dt}{\int_0^\iy w_{\al,\be,c;q}(s)\,ds}\,,
\quad{\rm where}\quad
w_{\al,\be,c;q}(t):=t^{c-1}\,
\frac{(-tq^{\be+1},-t^{-1}q^{\al+2};q)_\iy}{(-tq^{-c},-t^{-1}q^{1+c};q)_\iy}
\,.
\label{eq:3.07}
\eeq
Observe from \eqref{eq:3.01} that, for $k,l\in\Zplus$,
\begin{eqnarray*}
&&
\int_0^\iy \frac{t^k}{(-tq^{\be+1};q)_k}\,
\frac{t^{-l}}{(-t^{-1}q^{\al+2};q)_l}\,w_{\al,\be,c;q}(t)\,dt
=q^{-(k-l)(c+k-l)}
\int_0^\iy w_{\al+k,\be+l,c+k-l;q}(t)\,dt\qquad
\\
&&\qquad\qquad\qquad\qquad=
\frac{q^{c(c+k-l)}\,\Ga(c+k-l)\,\Ga(1-c-k+l)\,\Ga_q(\al+k+1)\,\Ga_q(\be+l+1)}
{\Ga_q(c+k-l)\,\Ga_q(1-c-k+l)\,\Ga_q(\al+\be+k+l+2)}
\\
&&\qquad\qquad\qquad\qquad=\frac{(q^{\al+1};q)_k\,(q^{\be+1};q)_l}
{q^{\half(k-l)(k-l-1)}\,(q^{\al+\be+2};q)_{k+l}}\,
\frac{q^{c^2}\,\Ga(c)\,\Ga(1-c)\,\Ga_q(\al+1)\,\Ga_q(\be+1)}
{\Ga_q(c)\,\Ga_q(1-c)\,\Ga_q(\al+\be+2)}\,.
\end{eqnarray*}
Thus
\beq
\int_0^\iy
\frac{q^{\half k(k-1)}\,t^k}
{(q^{\al+1},-tq^{\be+1};q)_k}\,
\frac{q^{\half l(l-1)}\,(t^{-1}q)^l}
{(q^{\be+1},-t^{-1}q^{\al+2};q)_l}\,
d\mu_{\al,\be,c;q}(t)=
\frac{q^{kl}}{(q^{\al+\be+2};q)_{k+l}}\,.
\label{eq:3.08}
\eeq
\bLP
{\bf Proof of Theorem \ref{th:1.13}}
\sLP
Multiply both sides of \eqref{eq:3.08} with
$\dstyle\frac{(q^{-n},q^{n+\al+\be+1};q)_k\,q^k}{(q;q)_k}$
and sum from $k=0$ to $n$.
Then the \RHS\ becomes
\[
\frac{1}{(q^{\al+\be+2};q)_l}\,
\qhyp21{q^{-n},q^{n+\al+\be+1}}{q^{\al+\be+l+2}}{q,q^{l+1}}
=\frac{(q^{l-n+1};q)_n}{(q^{\al+\be+2};q)_{l+n}}
=\frac{(q;q)_n\,\de_{n,l}}{(q^{\al+\be+2};q)_{2n}}\quad
(l=0,1,\ldots,n)
\]
by \cite[(1.5.2)]{GRah90}.
Thus
\bea
&&\int_0^\iy\qhyp22{q^{-n},q^{n+\al+\be+1}}
{q^{\al+1},-tq^{\be+1}}{q,-tq}
\frac{q^{l(l-1)/2}\,(tq^{-1})^l}
{(q^{\be+1},-t^{-1}q^{\al+2};q)_l}\,
d\mu_{\al,\be,c;q}(t)=
\frac{(q;q)_n\,\de_{n,l}}{(q^{\al+\be+2};q)_{2n}}\,,\qquad
\nonu
&&\qquad\qquad\qquad\qquad\qquad\qquad\qquad\qquad
\qquad\qquad\qquad\qquad\qquad\qquad l=0,1,\ldots,n,
\label{eq:3.09}
\\
&&\int_0^\iy
\frac{q^{k(k-1)/2}\,t^k}
{(q^{\al+1},-tq^{\be+1};q)_k}\,
\qhyp22{q^{-m},q^{m+\al+\be+1}}{q^{\be+1},-t^{-1}q^{\al+2}}
{q,-t^{-1}q^2} d\mu_{\al,\be,c;q}(t)=
\frac{(q;q)_m\,\de_{m,k}}{(q^{\al+\be+2};q)_{2m}}\,,
\nonu
&&\qquad\qquad\qquad\qquad\qquad\qquad\qquad\qquad
\qquad\qquad\qquad\qquad\qquad\qquad k=0,1,\ldots, m,
\label{eq:3.10}
\eea
where \eqref{eq:3.10} is obtained by a similar argument as \eqref{eq:3.09}.
Then \eqref{eq:3.09} and \eqref{eq:3.10} together with
\eqref{eq:1.14} imply the biorthogonality relations
\eqref{eq:1.16}.\qed
\bPP
The $q$-integral version of the Askey-Roy
integral \eqref{eq:3.01} is
\beq
\frac{(-sq^{-c},-s^{-1}q^{1+c};q)_\iy}{s^c\,(-s,-s^{-1}q;q)_\iy}\,
\int_0^{s\cdot\iy} t^{c-1}\,
\frac{(-tq^{\be+1},-t^{-1}q^{\al+2};q)_\iy}
{(-tq^{-c},-t^{-1}q^{1+c};q)_\iy}\,d_qt
=\frac{\Ga_q(\al+1)\,\Ga_q(\be+1)}
{\Ga_q(\al+\be+2)}\,.
\label{eq:3.11}
\eeq
Here $s>0$, $\Re\al,\Re\be>-1$, $c\in\RR$, and
the $q$-integral is defined by \eqref{eq:1.18}.
The case $s=q^c$ of \eqref{eq:3.11} (which is no real restriction)
was given in
\cite{G87} (see also \cite[(2.27)]{A88}) and in
\cite[Exercise 6.17(i)]{GRah90}).
Another approach to \eqref{eq:3.11} is presented in
\cite{DKa03}.

For $f$ any function on $(0,\iy)$ for which the sum below converges
absolutely, we have:
\bea
&&\frac{(-sq^{-c},-s^{-1}q^{1+c};q)_\iy}{s^c\,(-s,-s^{-1}q;q)_\iy}\,
\int_0^{s\cdot\iy} f(t)\,t^{c-1}\,
\frac{(-tq^{\be+1},-t^{-1}q^{\al+2};q)_\iy}
{(-tq^{-c},-t^{-1}q^{1+c};q)_\iy}\,d_qt
\nonu
&&\qquad=
\frac{(1-q)\,(-sq^{\be+1},-s^{-1}q^{\al+2};q)_\iy}
{(-s,-s^{-1}q;q)_\iy}\,
\sum_{k=-\iy}^\iy f(sq^k)\,
\frac{(-sq^{-\al-1};q)_k}{(-sq^{\be+1};q)_k}\,
q^{k(\al+1)}.
\label{eq:3.15}
\eea
The \RHS,
and thus the \LHS\ of \eqref{eq:3.15} is independent of $c$.
Henceforth we will take $c=0$ without loss of information.
Then \eqref{eq:3.11} together with \eqref{eq:3.15} takes the form
\bea
&&\int_0^{s\cdot\iy} t^{-1}\,
\frac{(-tq^{\be+1},-t^{-1}q^{\al+2};q)_\iy}
{(-t,-t^{-1}q;q)_\iy}\,d_qt
=\frac{(1-q)\,(-sq^{\be+1},-s^{-1}q^{\al+2};q)_\iy}
{(-s,-s^{-1}q;q)_\iy}
\nonu
&&\qquad\qquad\qquad\qquad
\times\sum_{k=-\iy}^\iy
\frac{(-sq^{-\al-1};q)_k}{(-sq^{\be+1};q)_k}\,
q^{k(\al+1)}
=\frac{\Ga_q(\al+1)\,\Ga_q(\be+1)}
{\Ga_q(\al+\be+2)}\,.
\label{eq:3.13}
\eea
The second equality in \eqref{eq:3.13} is Ramanujan's
${}_1\psi_1$ sum \cite[(5.2.1)]{GRah90}. This observation is the usual way
to prove \eqref{eq:3.11}.
With a completely analogous argument as used for the proof of Theorem
\ref{th:1.13}, we can next prove Theorem \ref{th:1.17}.
Details are omitted.
\mPP
In completion of this section,
observe the following symmetry of the functions
$\prat n\al\be t$:
\bea
&&\hskip-1.8cm\prat n\al\be t=
\qhyp22{q^{-n},q^{n+\al+\be+1}}{q^{\al+1},-tq^{\be+1}}{q,-tq}
\nonu
&&\qquad=\frac1{(-tq^{\be+1};q)_n}\,
\qhyp21{q^{-n},q^{-n-\be}}{q^{\al+1}}{q,-tq^{\be+n+1}}
\nonu\noalign{\allowbreak}
&&\qquad=\frac{(-1)^n\,(q^{\be+1};q)_n\,t^n}{(q^{\al+1},-tq^{\be+1};q)_n}\,
\qhyp21{q^{-n},q^{-n-\al}}{q^{\be+1}}{q,-t^{-1}q^{\al+n+1}}
\nonu\noalign{\allowbreak}
&&\qquad=\frac{(-1)^n\,(q^{\be+1};q)_n}{(q^{\al+1};q)_n}\,
\frac{t^n\,(-t^{-1}q^{\al+1};q)_n}{(-tq^{\be+1};q)_n}\,
\qhyp22{q^{-n},q^{n+\al+\be+1}}{q^{\be+1},-t^{-1}q^{\al+1}}{q,-t^{-1}q}
\nonu\noalign{\allowbreak}
&&\qquad=\frac{(-1)^n\,(q^{\be+1};q)_n}{(q^{\al+1};q)_n}\,
\frac{t^n\,(-t^{-1}q^{\al+1};q)_n}{(-tq^{\be+1};q)_n}\,
\prat n\be\al{t^{-1}}.
\label{eq:3.14}
\eea
Here the first and last equality are just \eqref{eq:1.14}. We use
\cite[(1.5.4)]{GRah90} for the second and fourth equality,
while the third equality is obtained by reversion of summation order in
a terminating $q$-hypergeometric series
(see \cite[Exercise 1.4.(ii)]{GRah90}).
\section{Concluding remarks} \label{sec:4}
%until 4.06
%
The results of this paper lead to several interesting questions.
I formulate two of these questions here. I also discuss specializations
of Theorem \ref{th:1.13}.
\PP
As I already mentioned in \S1, the new addition formula in the
case of the continuous $q$-Legendre polynomials
(formula \eqref{eq:1.06}) was first obtained
in a quantum group context, where a two-parameter family of
Askey-Wilson polynomials, including the continuous
$q$-Legendre polynomials, was interpreted
as spherical functions on the $SU_q(2)$ quantum group.
Here the left and right invariance of the spherical functions was
no longer \wrt the diagonal quantum subgroup, but infinitesimally \wrt
twisted primitive elements in the corresponding quantized universal
enveloping algebra. The $\si$ and $\tau$ variables in the addition
formula \eqref{eq:1.06} are parameters for the twisted primitive elements
occurring in the left respectively right invariance property.
On the other hand, the ${}_3\phi_2$ factors on the \RHS\ of
\eqref{eq:1.06}, involving $\si$ resp.\ $\tau$, can be rewritten
as the functions $\prat{n-k}kk{\,.\,}$ of argument $-q^{-\si}$
resp.\ $-q^{-\tau}$. So I wonder whether an
interpretation of these last functions and of their biorthogonality
(discussed in \S\ref{sec:3})  can be given
in the context of $SU_q(2)$.
\PP
A second question is whether the biorthogonality relations
for the functions $\prat n\al\be t$ (Theorems
\ref{th:1.13} and \ref{th:1.17}) fit into a more general class
of biorthogonal rational functions.
In fact, several papers have appeared during the last 10 or 20 years
which discuss explicit systems of biorthogonal rational functions
depending on many parameters and expressed as $q$-hypergeometric
functions, see
Rahman \cite{Rah86}, Wilson \cite{W91}, Ismail \& Masson \cite{IMas95} and
Spiridonov \& Zhedanov \cite{SpZhe00}.
However, I did not see how the functions $\prat n\al\be t$ and their
biorthogonality relations can be obtained as special or limit case
of families discussed in these references. If the functions
$\prat n\al\be t$ are indeed unrelated to the functions discussed
in these references, then it is a natural question how to generalize
the system of biorthogonal functions $\prat n\al\be t$.
\PP
On the other hand, our functions $\prat n\al\be t$ have some
interesting limit cases, one of which has occurred earlier
in literature. When we take limits for $\al\to\iy$ and/or
$\be\to\iy$ in \eqref{eq:1.14}, then we obtain:
\bea
\prat n\al\iy t&:=&{}_1\phi_1(q^{-n};q^{\al+1};q,-qt),
\label{eq:4.01}
\\\noalign{\allowbreak}
\prat n\iy\be t&:=&{}_1\phi_1(q^{-n};-tq^{\be+1};q,-qt),
\label{eq:4.02}
\\\noalign{\allowbreak}
\prat n\iy\iy t&:=&{}_1\phi_1(q^{-n};0;q,-qt).
\label{eq:4.03}
\eea
In \eqref{eq:4.01} and \eqref{eq:4.03} we have polynomials of degree $n$
in $t$,
rather than rational functions in $t$. The limit case $\be\to\iy$ of
the biorthogonality relations \eqref{eq:1.16} then becomes:
\beq
\int_0^\iy
\prat n\al\iy t\,\prat m\iy\al{t^{-1}q}\,
d\mu_{\al,\iy,c;q}(t)=
(-1)^n\,q^{-\half n(n-1)}\,(q;q)_n\,\de_{n,m},
\label{eq:4.04}
\eeq
where
\beq
d\mu_{\al,\iy,c;q}(t):=
\frac{q^{-c^2}\,\Ga_q(c)\,\Ga_q(1-c)\,(q^{\al+1};q)_\iy}
{\Ga(c)\,\Ga(1-c)\,(1-q)\,(q;q)_\iy}\,
\frac{t^{c-1}\,(-t^{-1}q^{\al+2};q)_\iy}
{(-tq^{-c},-t^{-1}q^{1+c};q)_\iy}\,dt\qquad(\al>-1).
\label{eq:4.05}
\eeq
The further limit case $\al\to\iy$ of
the biorthogonality relations \eqref{eq:4.04} then becomes
\bea
&&\hskip-1cm\frac{q^{-c^2}\,\Ga_q(c)\,\Ga_q(1-c)}
{\Ga(c)\,\Ga(1-c)\,(1-q)\,(q;q)_\iy}\,
\int_0^\iy
\prat n\iy\iy t\,\prat m\iy\iy{t^{-1}q}\,
\frac{t^{c-1}\,dt}
{(-tq^{-c},-t^{-1}q^{1+c};q)_\iy}
\nonu
&&\qquad\qquad\qquad\qquad\qquad\qquad\qquad\qquad
=(-1)^n\,q^{-\half n(n-1)}\,(q;q)_n\,\de_{n,m}.
\label{eq:4.06}
\eea
Similar limit cases can be considered for the biorthogonality relations
\eqref{eq:1.19}.
The biorthogonality relations \eqref{eq:4.06}
are essentially the ones
observed by Pastro \cite[pp.~532,~533]{P85}. He also points out that
biorthogonality relations of the form $\int P_n(t) Q_m(t^{-1})\,d\mu(t)=
h_n\,\de_{n,m}$
with $P_n$ and $Q_n$ polynomials of degree $n$ can be rewritten
as orthogonality relations on the real line for Laurent polynomials.
This is indeed the case in \eqref{eq:4.06}. Pastro also observes that
the biorthogonality measure occurring in \eqref{eq:4.06} is a
(non-unique) orthogonality measure for the Stieltjes-Wigert polynomials
(see \cite[\S3.27]{{KoekSw98}}).

\quad\\
\begin{footnotesize}
\begin{quote}
{ T. H. Koornwinder, Korteweg-de Vries Institute, University of
 Amsterdam,\\
 Plantage Muidergracht 24, 1018 TV Amsterdam, The Netherlands;

\vspace{\smallskipamount}
email: }{\tt thk@science.uva.nl}
\end{quote}
\end{footnotesize}
\end{document}